\documentclass[]{book}

\usepackage[utf8]{inputenc}
\usepackage{graphicx}
\usepackage{natbib}
\usepackage{xcolor}
\usepackage{algorithmic,algorithm}

\newtheorem{theorem}{Theorem}[section]

\newtheorem{remark}[theorem]{Remark}

\newtheorem{problem}[theorem]{Problem}
\newtheorem{definition}[theorem]{Definition}

\newcommand{\real}{{\mathbb{R}}}
\newcommand{\realpositive}{\mathbb{R}_{>0}}
\newcommand{\realnonnegative}{\mathbb{R}_{\ge 0}}
\newcommand{\integernonnegative}{\mathbb{Z}_{\ge 0}}
\newcommand{\integerpositive}{\mathbb{Z}_{> 0}}

\newcommand{\GG}{{\mathcal{G}}}

\newcommand{\NN}{{\mathcal{N}}}
\newcommand{\douti}{d_i^{\operatorname{out}}}
\newcommand{\Nouti}{\mathcal{N}_i^{\operatorname{out}}}
\newcommand{\Noutmax}{\mathcal{N}_{\max}^{\operatorname{out}}}
\newcommand{\Nini}{\mathcal{N}_i^{\operatorname{in}}}
\newcommand{\Dout}{D^{\operatorname{out}}}
\newcommand{\Din}{D^{\operatorname{in}}}

\newcommand{\BB}{{\mathcal{B}}}

\newcommand{\KL}{{\mathcal{KL}}}

\newcommand{\ones}[1]{\mathbf{1}_{#1}}
\newcommand{\zeros}[1]{\mathbf{0}_{#1}}

\newcommand{\diag}[1]{\operatorname{diag} ( #1)}

\newcommand{\vertices}{V}
\newcommand{\Norm}[1]{\|#1\|}
\newcommand{\Sym}[1]{{#1}_{s}}

\newcommand{\control}{\color{red!95!black}no}
\newcommand{\comm}{\color{green!60!black}yes}
\newcommand{\timeonly}{\color{red!95!black}time}
\newcommand{\state}{\color{green!60!black}state}
\newcommand{\centralized}{\color{red!95!black}centralized}
\newcommand{\decentralized}{\color{green!60!black}decentralized}
\newcommand{\undirected}{\color{red!95!black}undirected}
\newcommand{\directed}{\color{green!60!black}directed}
\newcommand{\continuous}{\color{red!95!black}continuous}
\newcommand{\periodic}{\color{green!60!black}periodic}
\newcommand{\zeno}{\color{red!95!black}no}
\newcommand{\nozeno}{\color{green!60!black}yes}

\newcommand{\until}[1]{\{1,\dots, #1\}}
\newcommand{\setdef}[2]{\{#1 \; | \; #2\}}
\newcommand{\TwoNorm}[1]{\|#1\|_2}
\renewcommand{\TwoNorm}[1]{\|#1\|}

\newcommand{\commgraph}{\GG_\text{comm}}

\renewcommand{\commgraph}{\GG}

\newcommand{\oprocendsymbol}{\hbox{$\bullet$}}
\newcommand{\oprocend}{\relax\ifmmode\else\unskip\hfill\fi\oprocendsymbol}

\usepackage{wileySTM}
\usepackage[\ifnum\pdfoutput=1breaklinks\fi]{hyperref}
\author{C. Nowzari, J. Cort\'es, and G. J. Pappas}
\title{Event-triggered communication and control for multi-agent average consensus}

\begin{document}
\frontmatter
\mainmatter

\section{Introduction}\label{se:intro}

In this chapter we look at one of the canonical driving examples for
multi-agent systems: average consensus. In this scenario, a group of
agents seek to agree on the average of their initial states. Depending
on the particular application, such states might correspond to sensor
measurements, estimates about the position of a target, or some other
data that needs to be fused.  Due to its numerous applications in
networked systems, many algorithmic solutions exist to the multi-agent
average consensus problem; however, a majority of them rely on agents
having continuous or periodic availability of information from other
agents.  Unfortunately, this assumption leads to inefficient
implementations in terms of energy consumption, communication
bandwidth, network congestion, and processor usage.  Motivated by
these observations, our main goal here is the design of provably
correct distributed event-triggered strategies that autonomously decide when
communication and control updates should occur so that the resulting
asynchronous network executions still achieve average consensus.

The literature and motivation behind multi-agent average consensus is
extensive, see e.g.,~[1-4] and
references therein. This chapter aims to provide a conceptual
introduction to event-triggered control strategies applied to
consensus problems. Triggered controllers seek to understand the
trade-offs between computation, communication, sensing, and actuator
effort in achieving a desired task with a guaranteed level of
performance.  Early works~[5] only consider tuning
controller executions to the state evolution of a given system, but
these ideas have since been extended to consider other
tasks such as when to take the sample of a state or when to broadcast
information over a wireless network;
see~[6] and references therein for a recent
overview.  Among the many references in the context of multi-agent
systems,~[7] specifies the responsibility of each agent
in updating the control signals,~[8] considers network
scenarios with disturbances, communication delays, and packet drops,
and~[9] studies decentralized event-based control that
incorporates estimators of the interconnection signals among agents.
These works are all concerned with designing event-triggers that ultimately determine
when control signals should be updated in addition to how. 
Several works have explored the application of event-triggered ideas
to the acquisition of information by the agents rather than only for actuation.
To this end,~[10-12] combine
event-triggered controller updates with sampled data that allows for
the periodic evaluation of the triggers. Instead, some works~[13] drop the
need for periodic access to information by considering event-based
broadcasts, where agents decide with local information only when to
obtain further information about neighbors. Self-triggered
control~[14,15] relaxes the need for local
information by deciding when a future sample of the state should be
taken based on the available information from the last sampled
state. Team-triggered coordination~[16] combines the
strengths of event- and self-triggered control into a unified approach
for networked systems.

\begin{table}
\begin{center}
\begin{tabular}{|l|c|c|c|c|c|c|}
\hline
& Triggered & Trigger  & Memory & Graph & Trigger & Provably \\
& comm? & dependence & structure & structure & evaluation & no Zeno?\\
\hline
[17] & \control & \state & \centralized & \undirected & \continuous & \nozeno \\ 
\hline
[17] & \control & \state & \decentralized & \undirected & \continuous & \zeno \\ 
\hline
[18] & \control & \state & \centralized & \directed & \continuous & \nozeno \\
\hline
[18] & \control & \state & \decentralized & \directed & \continuous & \zeno \\
\hline
[19] & \comm & \timeonly & \decentralized & \undirected & \continuous & \zeno \\
\hline
[19] & \comm & \timeonly & {\color{yellow!70!black}requires $\lambda_2$} & \undirected & \continuous & \nozeno \\
\hline
[20] & \comm & \state & {\color{green!60!black}requires $N$} & \undirected & \continuous & \zeno \\
\hline
[21] & \comm & \state & {\color{yellow!70!black}requires $K$} & \directed & \periodic & \nozeno \\
\hline
[12] & \comm & \state & \decentralized & \undirected & \periodic & \nozeno \\
\hline
[22,23] & \comm & \state & \decentralized & \directed & \periodic & \nozeno \\
\hline
\end{tabular}
\end{center}
\caption{Event-trigged multi-agent average consensus \vspace*{-5ex}}\label{ta:survey}
\end{table}

\paragraph{Organization}
Table~\ref{ta:survey} shows the progression of event-triggered
consensus problems that are covered in this chapter. It should be
noted that this is a very narrow scope on the field of event-triggered
consensus problems intended to introduce the high-level ideas behind
event-triggered communication and control laws and provide insight
into how they are designed.  In particular, this chapter only
discusses works that consider single-integrator dynamics and no
uncertainties (e.g., disturbances, noise, quantization, wireless
communication issues). Given that this is currently an active area of
research, it goes without saying that there are many important related
works that are not highlighted here. Examples include scenarios with
disturbances, sensor noise, delayed communication, quantized
communication, packet drops, more general dynamics, dynamic
topologies, and heterogeneous agents; to name a few. Lastly, it should
also be noted that although the table references journal articles that
first present these ideas going back to 2012, preliminary results from
these works has been presented at various conferences as early as
2008.  The contents of the chapter are summarized next.

The first application of event-triggered control to the multi-agent
consensus problem was in~[17], where the authors
propose a Lyapunov-based event-triggered control strategy that
dictates when agents should update their control signals.
Unfortunately, its implementation relies on each agent having perfect
information about their neighbors at all times. Identifying this
limitation, the authors in~[19] propose an
event-triggered communication and control law that not only determines
when agents should update their control signals, but also when they
should broadcast state information to their neighbors. However, the
drawback of the proposed algorithm is that it is a time-dependent
triggering rule with design parameters that are difficult to choose to
yield good performance. Instead, a state-dependent triggering rule is
proposed in~[20] which better aligns the events
with the desired task; this is explained in more detail later.
Lastly, all the above algorithms assume continuous evaluation of some
function is possible to determine exactly when some event has
occurred. Even in scenarios where Zeno behavior (an infinite number of
events occurring in a finite period of time) can be provably avoided,
the time between events may still be arbitrarily small which can be
problematic for digital implementations.  Consequently, the
works~[12,21-23]
propose algorithms that only require triggering functions to be
evaluated periodically rather than continuously. Finally, we close the
chapter by identifying some shortcomings of the current state of the
art and ideas for future work.

\section{Preliminaries}

This section introduces some notational conventions and notions on
graph theory.  Let $\real$, $\realpositive$, $\realnonnegative$, and
$\integerpositive$ denote the set of real, positive real, nonnegative
real, and positive integer numbers, respectively. We denote by
$\ones{N}$ and $\zeros{N} \in \real^N$ the column vectors with entries
all equal to one and zero, respectively. We let $\Norm{\cdot}$ denote
the Euclidean norm on $\real^N$. We let $\diag{\real^N} = \setdef{x
  \in \real^N}{x_1=\dots=x_N} \subset \real^N$ be the agreement
subspace in $\real^N$.  For a finite set $S$, we let $|S|$ denote its
cardinality.
Given $x, y \in \real$, Young's inequality
states that, for any $\epsilon \in \realpositive$,
\begin{align}\label{eq:Young}
  xy \leq \frac{x^2}{2\epsilon} + \frac{ \epsilon y^2}{2}.
\end{align}
A weighted directed graph (or weighted digraph) $\commgraph = (V,E,W)$
is comprised of a set of vertices $V = \until{N}$, directed edges $E
\subset V \times V$ and weighted adjacency matrix $W \in
\realnonnegative^{N \times N}$. Given an edge $(i,j) \in E$, we refer
to $j$ as an out-neighbor of $i$ and $i$ as an in-neighbor of $j$.
The sets of out- and in-neighbors of a given node $i$ are $\Nouti$ and
$\Nini$, respectively. The weighted adjacency matrix $W \in \real^{N
  \times N}$ satisfies $w_{ij} > 0$ if $(i, j) \in E$ and $w_{ij} = 0$
otherwise. The graph~$\commgraph$ is \emph{undirected} if and only
if $w_{ij} = w_{ji}$ for all $i,j \in \vertices$.
A path from vertex $i$ to $j$ is an ordered sequence of
vertices such that each intermediate pair of vertices is an edge.  A
digraph $\commgraph$ is strongly connected if there exists a path from
all $i \in V$ to all $j \in V$.  The out- and in-degree matrices
$\Dout$ and $\Din$ are diagonal matrices where
\begin{align*}
  \douti = \sum_{j \in \Nouti} w_{ij} , \quad d_i^{\operatorname{in}} =
  \sum_{j \in \Nini} w_{ji} ,
\end{align*}
respectively. A digraph is weight-balanced if $\Dout = \Din$.  The
(weighted) Laplacian matrix is $L = \Dout - W$. Based on the structure
of $L$, at least one of its eigenvalues is zero and the rest of them
have nonnegative real parts.  If the digraph $\commgraph$ is strongly
connected, $0$ is a simple eigenvalue with associated eigenvector
$\ones{N}$.  The digraph~$\commgraph$ is weight-balanced if and only
if $\ones{N}^T L = \zeros{N}$ if and only if $ \Sym{L}=\frac{1}{2}(L +
L^T)$ is positive semidefinite.  For a strongly connected and
weight-balanced digraph, zero is a simple eigenvalue of $\Sym{L} $.
In this case, we order its eigenvalues as $\lambda_1 = 0<\lambda_2\leq
\dots \leq \lambda_N$, and note the inequality
\begin{align}\label{eq:LapBound}
  x^T L x \ge \lambda_2(\Sym{L}) \| x -\frac{1}{N} (\ones{N}^T x)
  \ones{N} \|^2 ,
\end{align}
for all $x \in \real^N$.  The following property will also be of use
later,
\begin{align}\label{eq:lap2-bound}
  \lambda_2(\Sym{L}) x^T L x \leq x^T \Sym{L}^2 x \leq
  \lambda_N(\Sym{L}) x^T L x .
\end{align}
This can be seen by noting that $\Sym{L}$ is diagonalizable and
rewriting $\Sym{L} = S^{-1} D S$, where $D$ is a diagonal matrix
containing the eigenvalues of~$\Sym{L}$.

%

\subsection{Event-triggered control of linear systems}\label{se:event}

Here we provide a very basic working introduction to the general idea of event-triggered control
by working through a simple linear control problem. The exposition closely
follows~[24]. The remainder of this chapter then focuses
on how these elementary ideas are extended to be applied to much more in the context of multi-agent
consensus on networks. We refer the interested reader to~[6] for further details on the subject of event-triggered control in general.

Consider a linear control system
\begin{align}\label{eq:lindynamics}
  \dot{x} = Ax + Bu,
\end{align}
with $ x \in \real^n$ and $u \in \real^m$. Our starting point
is the availability of a linear feedback controller $u^* = Kx$ such that the closed-loop system
\begin{align*}
  \dot{x} = (A+BK)x ,
\end{align*}
is asymptotically stable. Given a positive definite matrix $Q \in
\real^{n \times n}$, let $P \in \real^{n \times n}$ be the unique
solution to the Lyapunov equation $ (A+BK)^T P + P(A+BK) = -Q$. Then,
the evolution of the Lyapunov function $ V_c(x) = x^T P x$
along the trajectories of the closed-loop system is
\begin{align*}
  \dot{V}_c &= x^T ( (A+BK)^T P + P(A+BK) ) x = -x^T Q x.
\end{align*}

Consider now a sample-and-hold implementation of the controller, where
the input is not updated continuously, but instead at a sequence of
to-be-determined times $\{ t_\ell \}_{\ell \in \integernonnegative}
\subset \realnonnegative$,
\begin{align}\label{eq:zohcontroller}
  u(t) = Kx(t_\ell) , \quad t \in [t_\ell, t_{\ell+1} ).
\end{align}
Such an implementation makes sense in practical scenarios given the
inherent nature of digital systems. With this controller
implementation, the closed-loop system can be written as
\begin{align*}
  \dot{x} = (A+BK)x + BKe,
\end{align*}
where $ e(t) = x(t_\ell) - x(t)$, $ t \in [t_\ell, t_{\ell+1} )$, is
the state error. Then, the objective is to determine the sequence 
of times $\{t_\ell \}_{\ell \in \integernonnegative}$ to 
guarantee some desired level of performance for the
resulting system. To make this concrete, define the function
\begin{align*}
  V(t,x_0) = x(t)^T P x(t) ,
\end{align*}
for a given initial condition $x(0) = x_0$ (here, $t\mapsto x(t)$
denotes the evolution of the closed-loop system
using~\eqref{eq:zohcontroller}). We define the performance of the
system via a function $S : \realnonnegative \times \real^n
\rightarrow \realnonnegative$ that upper bounds the evolution
of~$V$. Then, the sequence of times $\{ t_\ell \}$ can be
implicitly defined as the times at which
\begin{align}\label{eq:trigger}
  V(t,x_0) \leq S(t,x_0)
\end{align}
is not satisfied. More specifically, this is an \emph{event-triggered} condition
that updates the actuator signal whenever $V(t_\ell, x_0) = S(t_\ell,
x_0)$. Assuming solutions are well defined, it is not difficult to
see that if the performance function satisfies $S(t,x_0) \leq
\beta(t,|x_0|)$, for some $\beta \in \KL$, then the closed-loop system
is globally uniformly asymptotically stable. Moreover, if $\beta$ is
an exponential function, the system is globally uniformly
exponentially stable.

Therefore, one only needs to guarantee the lack of Zeno behavior. We
do this by choosing the performance function $S$ so that the
inter-event times $t_{\ell+1} - t_\ell$ are lower bounded by some
constant positive quantity. This can be done in a number of ways.  For
the linear system~\eqref{eq:lindynamics}, it turns out that it is
sufficient to select $S$ satisfying $\dot{V}(t_\ell) <
\dot{S}(t_\ell)$ at the event times $t_\ell$ (this fact is formally
stated below in Theorem~\ref{th:lower-bound}).  To do so, choose $R
\in \real^{n \times n}$ positive definite such that $Q - R $ is also
positive definite. Then, there exists a Hurwitz matrix $A_s \in
\real^{n \times n}$ such that the Lyapunov equation
\begin{align*}
  A_s^T P + P A_s = -R 
\end{align*}
holds.  Consider the hybrid system,
\begin{align*}
  \dot{x}_s &= A_s x_s, \quad t \in [t_\ell, t_{\ell+1}),
  \\
  x_s(t_\ell) &= x(t_\ell),
\end{align*}
whose trajectories we denote by $t \mapsto x_s(t)$, and define the
performance function $S$ by
\begin{align*}
  S(t) = x_s^T(t) P x_s(t) .
\end{align*}

%
%
%

Letting $y = [x^T, e^T]^T \in \real^n \times \real^n$, we write the
continuous-time dynamics as
\begin{align*}
  \dot{y} = Fy , \quad t \in [t_\ell,t_{\ell+1}) ,
\end{align*}
where
\begin{align*}
  F = \left[
    \begin{array}{cc}
      A+BK & BK
      \\
      -A - BK & -BK
    \end{array}
  \right] .
\end{align*}
With a slight abuse of notation, we let $y_\ell = [x^T(t_\ell),
0^T]^T$ be the state $y$ at time $t_\ell$. Note that $e(t_\ell) = 0$,
for all $\ell \in \integernonnegative$, by definition of the update
times. With this notation, we can rewrite
\begin{align*}
  S(t) &= (C e^{F_s (t-t_\ell)} y_\ell)^T P (C e^{F_s (t-t_\ell)}
  y_\ell),
  \\
  V(t) &= (C e^{F (t-t_\ell)} y_\ell)^T P (C e^{F (t-t_\ell)} y_\ell)
  ,
\end{align*}
where
\begin{align*}
  F_s =
  \begin{bmatrix}
          A_s & 0
      \\
      0 & 0
  \end{bmatrix}
  , \quad C =
  \begin{bmatrix}
    I & 0
  \end{bmatrix}
.
\end{align*}
The condition~\eqref{eq:trigger} can then be rewritten as
\begin{align*}
  f(t,y_\ell) = y_\ell^T ( e^{F^T (t-t_\ell)} C^T P C e^{F (t-t_\ell)}
  - e^{F_s^T (t-t_\ell)} C^T P C e^{F_s (t-t_\ell)} )y_\ell \leq 0 .
\end{align*}
Note that because we consider a deterministic system here, with the 
information available at time $t_\ell$,
it is possible to determine the next time $t_{\ell+1}$ at
which~\eqref{eq:trigger} is violated by computing $t_{\ell+1} = h(x(t_\ell))$
as the time for which
\begin{align}\label{eq:selftrigger}
  f(h(x(t_\ell)), y_\ell) = 0.
\end{align}
The following result from~[24] provides a uniform
lower bound $t_\text{min}$ on the inter-event times $\{t_{\ell+1} -
t_\ell\}_{\ell \in \integernonnegative}$.

\begin{theorem}[Lower bound on inter-event times for event-triggered
  approach]\label{th:lower-bound}
  Given the system~\eqref{eq:lindynamics} with
  controller~\eqref{eq:zohcontroller} and controller updates given by
  the event-triggered policy~\eqref{eq:selftrigger}, the inter-event
  times are lower bounded by
  \begin{align*}
    t_{\operatorname{min}} = \min \setdef{t \in \realpositive}{
      \det(M(t)) = 0 } > 0,
  \end{align*}
  where
  \begin{align*}
    M(t) =
    \begin{bmatrix}
      I & 0
    \end{bmatrix}
    \left( e^{Ft}C^TPCe^{Ft}-e^{F_st}C^TPCe^{F_st} \right)
    \begin{bmatrix}
      I \\ 0
    \end{bmatrix}
    .
  \end{align*}
\end{theorem}
\smallskip

Note that the above result can also be interpreted in the context of a
periodic controller implementation: any period less than or equal to
$t_\text{min}$ results in a closed-loop system with asymptotic
stability guarantees.

\section{Problem statement}\label{se:statement}

%
%
%

We let $\commgraph$ denote the connected, undirected
communication graph that describes the communication topology in a network
of~$N$ agents. In other words, agent~$j$ can communicate with agent~$i$ if
$j$ is a neighbor of~$i$ in~$\commgraph$.  We denote by $x_i \in
\real$ the state of agent $i \in \until{N}$ and consider
single-integrator dynamics
\begin{align}\label{eq:dynamics}
\dot{x}_i(t) = u_i(t).
\end{align}
Then, the distributed controller
\begin{align}\label{eq:ideal}
  u_i^*(x) = -\sum_{j \in \NN_i} \left( x_i - x_j \right)
\end{align}
is known to drive the states of all agents to the average of the
initial conditions~[25,1]. This is
formalized in Theorem~\ref{th:main}.

\begin{theorem}[Continuous controller]\label{th:main}
  Given the dynamics~\eqref{eq:dynamics}, if all agents implement the
  control law~\eqref{eq:ideal}, then multi-agent average consensus is
  achieved; i.e.,
  \begin{align}\label{eq:averageconsensus}
    \lim_{t \rightarrow \infty} x_i(t) = \frac{1}{N} \sum_{j=1}^N
    x_j(0)
  \end{align}
  for all $i \in \until{N}$.
\end{theorem}

Unfortunately, implementing~\eqref{eq:ideal} in a digital setting is
not possible since it requires all agents to have continuous access to
the state of their neighbors and the control inputs $u_i(t)$ must also
be updated continuously. This is especially troublesome in the context
of wireless network systems since this means agents must communicate
with each other continuously as well. Instead, this chapter is
interested in event-triggered communication and control strategies to
relax these requirements.

\section{Centralized event-triggered control}

Consider the dynamics~\eqref{eq:dynamics} and the ideal control
law~\eqref{eq:ideal}.  Letting $x = (x_1, \dots, x_N)^T$ and $u =
(u_1, \dots, u_N)^T$, the closed-loop dynamics of the ideal system is
given by
\begin{align}\label{eq:idealfull}
  \dot{x}(t) = -Lx(t) .
\end{align}
As stated before, implementing this requires all agents to
continuously update their control signals which is not realistic for
digital controllers. Instead, following the basic idea for event-triggered
control presented in Section~\ref{se:event}, let us consider a digital implementation
of this ideal controller
\begin{align}\label{eq:centralcontrol}
  u(t) = -Lx(t_\ell), \quad t \in [t_\ell, t_{\ell+1} ),
\end{align}
where the event times $\{ t_\ell \}_{\ell \in \integernonnegative}$
are to be determined such that the system still converges to the
desired state. Let $e(t) = x(t_\ell) - x(t)$ for $t \in [t_\ell,
t_{\ell+1})$ be the state measurement error.  For simplicity, we
denote by $\hat{x}(t) = \hat{x}(t_\ell)$ for $t \in [t_\ell,
t_{\ell+1})$ as the state that was used in the last computation of the
control signal.  The closed-loop dynamics of the
controller~\eqref{eq:centralcontrol} is then given by
\begin{align}\label{eq:centraldynamics}
  \dot{x}(t) = -L \hat{x}(t) = -L( x(t) + e(t) ).
\end{align}

The problem can now be formalized as follows.

\begin{problem}[Centralized event-triggered control]
  Given the closed-loop dynamics~\eqref{eq:centraldynamics}, find an
  event-trigger such that the sequence of times~$\{ t_\ell \}_{\ell
    \in \integernonnegative}$ ensures multi-agent average
  consensus~\eqref{eq:averageconsensus} is achieved.
\end{problem}

Following~[17], to solve this problem we consider the
Lyapunov function
\begin{align*}
  V(x) = x^T L x.
\end{align*}
Given the closed-loop dynamics~\eqref{eq:centraldynamics}, we have
\begin{align*}
  \dot{V} = x^T L \dot{x} = -x^T L L (x + e) =
  -\underbrace{\TwoNorm{Lx}^2}_\text{"good"} - \underbrace{x^T L L
    e}_\text{"bad"} .
\end{align*}
For simplicity, we are not interested in characterizing any specific
performance as in Section~\ref{se:event}. Instead, we are only interested
in asymptotic stability. The main idea of event-triggered control is 
then to determine when the
controller should be updated (i.e., when $e$ should be set to 0) by
balancing the ``good'' term against the ``bad'' term. More specifically,
we are interested in finding conditions on the error~$e$ such that $\dot{V} <
0$ at all times. Using norms, we can bound
\begin{align*}
  \dot{V} \leq -\TwoNorm{Lx}^2 + \TwoNorm{Lx}\TwoNorm{L}\TwoNorm{e}.
\end{align*}
Then, if we enforce the error~$e$ to satisfy
\begin{align*}
  \TwoNorm{e} \leq \sigma \frac{ \TwoNorm{Lx} }{\TwoNorm{L}} ,
\end{align*}
with $\sigma \in (0,1)$ for all times, we have
\begin{align*}
  \dot{V} \leq (\sigma - 1) \TwoNorm{Lx}^2,
\end{align*}
which is strictly negative for all $Lx \neq 0$. The following
centralized event-trigger ensures this is satisfied at all times.

\begin{theorem}[Centralized event-triggered control]
  Given the closed-loop dynamics~\eqref{eq:centraldynamics}, if the
  update times are determined as the times when
  \begin{align}\label{eq:centraltrigger}
    f(x,e) \triangleq \TwoNorm{e} - \sigma \frac{
      \TwoNorm{Lx}}{\TwoNorm{L}} = 0,
  \end{align}
  then the system achieves multi-agent average consensus.
\end{theorem}

In other words, given a control update at time~$t_\ell$, the next time
$t_{\ell+1}$ is given by
\begin{align*}
  t_{\ell+1} = \min \setdef{t' > t_\ell}{\TwoNorm{e(t')} = \sigma
    \frac{\TwoNorm{Lx(t')}}{\TwoNorm{L}}}.
\end{align*}
The algorithm is formalized in Table~\ref{tab:algorithm2}.

\begin{table}[htb]
  \centering
  \framebox[.9\linewidth]{\parbox{.85\linewidth}{%
      \parbox{\linewidth}{At times $t \in [t_\ell,t_{\ell+1})$, system (continuously) performs:}
      \vspace*{-1.5ex}
      \begin{algorithmic}[1]
      \STATE set $\hat{x}(t) = x(t_\ell)$
      \STATE set $e(t) = \hat{x}(t) - x(t)$
        \IF{$\TwoNorm{e(t)} = \sigma \frac{\TwoNorm{L x(t)}}{\TwoNorm{L}}$}  
        \STATE set $t_{\ell+1} = t$
        \STATE set $\hat{x}(t) = x_i(t_{\ell+1})$
        \STATE set $\ell = \ell + 1$
        \ENDIF
        \STATE set $u(t) = - L \hat{x}(t)$
      \end{algorithmic}}}
  \caption{Centralized event-triggered control.}\label{tab:algorithm2}
\end{table}

The proof of convergence to the desired state then follows directly
from the proof of Theorem~\ref{th:main} and the fact that the sum of
all states is still an invariant quantity. Furthermore, the authors
in~[17] are able to rule out the existence of Zeno
behavior (formally defined below) by showing there exists a positive time
\begin{align*}
  \tau = \frac{\sigma}{\TwoNorm{L}(1+\sigma)}
\end{align*}
bounding the inter-event times, i.e.,
\begin{align*}
  t_{\ell+1} - t_\ell \geq \tau
\end{align*}
for all $\ell \in \integernonnegative$.

\begin{definition}[Zeno behavior]
    If there exists $T > 0$ such that $t_\ell
    \leq T$ for all $\ell \in \integernonnegative$, then the system is
    said to exhibit \emph{Zeno behavior}.
\end{definition}

The centralized event-triggered controller~\eqref{eq:centralcontrol}
with triggering law~\eqref{eq:centraltrigger} relaxes the requirement
that agents need to continuously update their control signals, but it
still has many issues. One of them is that the event-trigger~$f(x,e)$
requires full state information to implement. Next, we provide a
distributed solution instead of a centralized one.

\section{Decentralized event-triggered control}

In the previous section we presented a centralized event-triggered
control law to solve the multi-agent average consensus
problem. Unfortunately, implementing this requires a centralized
decision maker and requires all agents in the network to update their
control signals simultaneously. In this section we relax this
requirement by following~[17].

Let us now consider a distributed digital implementation of the ideal
controller~\eqref{eq:ideal}. In this case we assume each agent~$i$ has
its own sequence of event times $\{ t_\ell^i \}_{\ell \in
  \integernonnegative}$.  At any given time~$t$, let $\hat{x}_i(t) =
x_i(t_\ell^i)$ for $t \in [t_\ell^i, t_{\ell+1}^i )$ be the state of
agent~$i$ at its last update time. The distributed event-triggered
controller is then given by
\begin{align}\label{eq:decentralcontrol}
  u_i(t) = - \sum_{j \in \NN_i} (\hat{x}_i(t) - \hat{x}_j(t) ) .
\end{align}
It is important to note here that the latest updated state
$\hat{x}_j(t)$ of agent $j \in \NN_i$ appears in the control signal
for agent~$i$. This means that when an event is triggered by a
neighboring agent~$j$, agent~$i$ also updates its control signal
accordingly.  As in the centralized case, let $e_i(t) = x_i(t_\ell^i)
- x_i(t)$ be the state measurement error for agent~$i$. Then, letting
$\hat{x} = (\hat{x}_1, \dots, \hat{x}_N)^T$ and $e = (e_1, \dots,
e_N)^T$, the closed-loop dynamics of the
controller~\eqref{eq:decentralcontrol} is given by
\begin{align}\label{eq:decentraldynamics} 
  \dot{x}(t) = - L \hat{x}(t) = - L ( x(t) + e(t) ).
\end{align}
The problem can now be formalized as follows.

\begin{problem}[Decentralized event-triggered control]
  Given the closed-loop dynamics~\eqref{eq:decentraldynamics}, find an
  event-trigger for each agent~$i$ such that the sequence of times $\{
  t_\ell^i \}_{\ell \in \integernonnegative}$ ensures multi-agent
  average consensus~\eqref{eq:averageconsensus} is achieved.
\end{problem}

Following~[17], to solve this problem we again
consider the Lyapunov function
\begin{align*}
  V(x) = x^T L x.
\end{align*}
Given the closed-loop dynamics~\eqref{eq:decentraldynamics}, we have
\begin{align*}
  \dot{V} = -\TwoNorm{Lx}^2 - x^T L L e .
\end{align*}
As before, we are interested in finding conditions on the error~$e$ such
that~$\dot{V} < 0$ at all times; however, we must now do this in a
distributed way. For simplicity, let $Lx \triangleq z = (z_1, \dots,
z_N)^T$. Then, expanding out $\dot{V}$ yields
\begin{align*}
  \dot{V} &= - \sum_{i=1}^N z_i^2 - \sum_{j \in \NN_i} z_i (e_i - e_j) \\
  &= -\sum_{i=1}^N z_i^2 - | \NN_i | z_i e_i + \sum_{j \in \NN_i} z_i
  e_j .
\end{align*}
Using Young's inequality~\eqref{eq:Young} and the fact
that~$\commgraph$ is symmetric, we can bound this by
\begin{align*}
  \dot{V} \leq - \sum_{i=1}^N (1 - a|\NN_i|)z_i^2 + \frac{1}{a}
  |\NN_i| e_i^2
\end{align*}
for all $a > 0$. Letting $a \in (0, 1/|\NN_i|)$ for all $i$, if we can
enforce the error of all agents to satisfy
\begin{align*}
  e_i^2 \leq \frac{ \sigma_i a (1 - a|\NN_i|)}{|\NN_i|} z_i^2
\end{align*}
with $\sigma_i \in (0,1)$ for all times, we have
\begin{align*}
  \dot{V} \leq \sum_{i=1}^N (\sigma_i - 1)(1-a|\NN_i|)z_i^2,
\end{align*}
which is strictly negative for all $Lx \neq 0$. The following
decentralized event-trigger ensures this is satisfied at all times.

\begin{theorem}[Decentralized event-triggered
  control]\label{th:event-control}
  Given the closed-loop dynamics~\eqref{eq:decentraldynamics}, if the
  updates times of each agent~$i$ are determined as the times when
  \begin{align}\label{eq:decentraltrigger}
    f_i(x_i, e_i, \{ x_j \}_{j \in \NN_i} ) \triangleq e_i^2 - \frac{
      \sigma_i a (1 - a|\NN_i|)}{|\NN_i|} z_i^2 = 0,
  \end{align}
  with $0 < a < 1/|\NN_i|$, then the system achieves multi-agent
  average consensus.
\end{theorem}

Note that the trigger~\eqref{eq:decentraltrigger} can be evaluated by
agent~$i$ using only information about its own and neighbors' states.
The algorithm is formalized in Table~\ref{tab:algorithm3}. 	

\begin{table}[htb]
  \centering
  \framebox[.9\linewidth]{\parbox{.85\linewidth}{%
      \parbox{\linewidth}{At times $t \in [t_\ell^i,t_{\ell+1}^i)$, agent~$i$ (continuously) performs:}
      \vspace*{-1.5ex}
      \begin{algorithmic}[1]
      \STATE set $z_i(t) = \sum_{j \in \NN_i} (x_i(t) - x_j(t))$
      \STATE set $e_i(t) = \hat{x}_i(t) - x_i(t)$
        \IF{$e_i(t)^2 = \frac{\sigma_i a (1 - a | \NN_i |)}{| \NN_i |} z_i(t)^2$}  
        \STATE set $t_{\ell+1}^i = t$
        \STATE broadcast $\hat{x}_i(t) = x_i(t_{\ell+1}^i)$ to neighbors $j \in \NN_i$
        \STATE set $\ell = \ell + 1$
        \ENDIF
      \STATE set $u_i(t) = - \sum_{j \in \NN_i} (\hat{x}_i(t) - \hat{x}_j(t))$
      \end{algorithmic}}}
  \caption{Decentralized event-triggered control.}\label{tab:algorithm3}
\end{table}

The proof of convergence to the desired state then directly follows
from the proof of Theorem~\ref{th:main} and the fact that the sum of
all states is still an invariant quantity. Furthermore, the authors
in~[17] are able to show that at all times there
exists one agent~$i$ for which the inter-event times are strictly
positive. Unfortunately, this is not enough to rule out Zeno behavior
which is quite problematic, both from a pragmatic and theoretical
viewpoint, as the trajectories of the system are no longer well-defined
beyond the accumulation point in time.

\begin{remark}[Convergence and Zeno behavior]
{\rm It should be noted here that when we refer to a ``proof of convergence''
for any closed-loop dynamics, it is only valid for trajectories do no exhibit
Zeno behavior. Consequently, being able to guarantee Zeno behaviors do not
occur is extremely important in validating the correctness of a given algorithm.
We formalize the definition of Zeno behavior next. } \oprocend
\end{remark}

\begin{remark}[Directed graphs] {\rm All the work
    from~[17] has also been extended to consider
    weight-balanced directed graphs in~[18]. For
    brevity, we defer the discussion on directed graphs to
    Section~\ref{se:directed}.  }
\end{remark}

The decentralized event-triggered
controller~\eqref{eq:decentralcontrol} with triggering
law~\eqref{eq:decentraltrigger} relaxes the requirement that agents
need to continuously update their control signals; however, there are
still some severe issues. Although each agent now has a local
event-triggering condition, it requires continuous information about
all of its neighbors to implement it. This is still troublesome in a
wireless network setting where this implies continuous communication
among agents is still required.  We address this next.

\section{Decentralized event-triggered communication and control}\label{se:coordination}

In the previous sections we presented event-triggered \emph{control}
laws to determine when control signals should be updated; however,
this relied on the continuous availability of some state
information. In particular, each agent~$i$ requires exact state
information about their neighbors~$j \in \NN_i$ to evaluate the
trigger~\eqref{eq:decentraltrigger} and determine when its control
signal $u_i$ should be updated. Instead, we are now interested in
developing event-triggered \emph{communication and control} laws such
that each agent~$i$ must not only determine when to update its control
signal but also when to communicate with its neighbors. For
simplicity, we refer to communication and control together as
`coordination.'

As in the previous section, we assume each agent~$i$ has its own
sequence of event times $\{ t_\ell^i \}_{\ell \in
  \integernonnegative}$.  However, these update times now correspond
to when messages are broadcast in addition to when control signals are
updated. At any given time~$t$, let $\hat{x}_i(t) = x_i(t_\ell^i)$ for
$t \in [t_\ell^i, t_{\ell+1}^i )$ be the last broadcast state of
agent~$i$. Then, at any given time~$t$, agent $i$ only has access to
the last broadcast state $\hat{x}_j(t)$ of its neighbors $j \in \NN_i$
rather than exact states $x_j(t)$.

The distributed event-triggered controller is then still given by
\begin{align}\label{eq:decentralcontrol2}
  u_i(t) = - \sum_{j \in \NN_i} (\hat{x}_i(t) - \hat{x}_j(t) ) .
\end{align}
It is important to note here that the latest broadcast state
$\hat{x}_j(t)$ of agent $j \in \NN_i$ appears in the control signal
for agent~$i$. This means that when an event is triggered by a
neighboring agent~$j$, agent~$i$ also updates its control signal
accordingly.  As before, let $e_i(t) = x_i(t_\ell^i) - x_i(t)$ be the
state measurement error for agent~$i$. Then, letting $\hat{x} =
(\hat{x}_1, \dots, \hat{x}_N)^T$ and $e = (e_1, \dots, e_N)^T$, the
closed-loop dynamics of the controller~\eqref{eq:decentralcontrol} is
given by
\begin{align}\label{eq:decentraldynamics2}
  \dot{x}(t) = - L \hat{x}(t) = - L ( x(t) + e(t) ).
\end{align}
The problem can now be formalized as follows. However, it should be
noted that we are now looking for a strictly \emph{local}
event-trigger for each agent~$i$ that doesn't require exact
information about its neighbors. More specifically, we recall the
result of Theorem~\ref{th:event-control} and notice that the
event-trigger for agent~$i$ depends on the exact state $x_j(t)$ of all
its neighbors $j \in \NN_i$. In this section we are interested in
finding a trigger that only depends on the last broadcast information
$\hat{x}_j(t)$ instead.

\begin{problem}[Decentralized event-triggered coordination]
  Given the closed-loop dynamics~\eqref{eq:decentraldynamics}, find a
  local event-trigger for each agent~$i$ such that the sequence of
  times $\{ t_\ell^i \}_{\ell \in \integernonnegative}$ ensures
  multi-agent average consensus~\eqref{eq:averageconsensus} is
  achieved.
\end{problem}

Here we present two classes of event-triggered coordination solutions
to the problem above: time-dependent and state-dependent triggers. The
time-dependent event-trigger to solve this problem was first developed
in~[19] and is presented next. The algorithm is
formalized in Table~\ref{tab:algorithm4}.

\begin{theorem}[Decentralized event-triggered coordination
  (time-dependent)]\label{th:time}
  Given the closed-loop dynamics~\eqref{eq:decentraldynamics}, if the
  updates times of each agent~$i$ are determined as the times when
  \begin{align}
    f_i(e_i(t),t) \triangleq \TwoNorm{e_i(t)} - (c_0 + c_1 e^{-\alpha
      t} ) = 0,
  \end{align}
  with constants $c_0, c_1 \geq 0$ and $c_0 + c_1 > 0$, then the
  system reaches a neighborhood of multi-agent average consensus
  upper-bounded by
  \begin{align*}
    r = \TwoNorm{L} \sqrt{N} c_0 / \lambda_2(L) .
  \end{align*}
  Moreover, if $c_0 > 0$ or $0 < \alpha < \lambda_2(L)$, then the closed-loop
  system does not exhibit Zeno behavior.
\end{theorem}

\begin{table}[htb]
  \centering
  \framebox[.9\linewidth]{\parbox{.85\linewidth}{%
      \parbox{\linewidth}{At times $t \in [t_\ell^i,t_{\ell+1}^i)$, agent~$i$ (continuously) performs:}
      \vspace*{-1.5ex}
      \begin{algorithmic}[1]
      \STATE set $\hat{x}_i(t) = x_i(t_\ell^i)$
      \STATE set $e_i(t) = \hat{x}_i(t) - x_i(t)$
        \IF{$|e_i(t)| = c_0 + c_1 e^{- \alpha t}$}  
        \STATE set $t_{\ell+1}^i = t$
        \STATE broadcast $\hat{x}_i(t) = x_i(t_{\ell+1}^i)$ to neighbors $j \in \NN_i$
        \STATE set $\ell = \ell + 1$
        \ENDIF
      \STATE set $u_i(t) = - \sum_{j \in \NN_i} (\hat{x}_i(t) - \hat{x}_j(t))$
      \end{algorithmic}}}
  \caption{Decentralized event-triggered coordination (time-dependent).}\label{tab:algorithm4}
\end{table}

The proof of convergence is shown in the appendix; however, we are now
also interested in guaranteeing Zeno behavior does not occur to verify
the correctness of the algorithm as mentioned earlier.

The main drawback of the event-triggered communication and control law
proposed in Theorem~\ref{th:time} is that although the parameters
$c_0, c_1,$ and $\alpha$ play very important roles in the performance
of the algorithm (e.g., convergence speed and amount of triggers),
there is no good way of choosing these parameters a priori, without
any global knowledge. Furthermore, the initial condition also plays an
important role in the performance of the algorithm.

In particular we focus our discussion here on the parameters $c_0$ and
$\alpha$ and their effects on convergence and possible Zeno
behaviors. We begin with the more desirable $c_0 = 0$ case, as in this
case the result of Theorem~\ref{th:time} states that the system will
asymptotically achieve exact multi-agent average consensus as defined
in~\eqref{eq:averageconsensus}. However, in this case we require
$\alpha < \lambda_2(L)$ to guarantee Zeno behaviors can be avoided
and, unfortunately, $\lambda_2(L)$ is a global quantity that requires
knowledge about the entire communication topology to
compute. {There are indeed methods for estimating
  this quantity in a distributed way (see e.g.,~[26,27]), but we do not discuss this
  here.}
On the other hand, when $c_0 > 0$ we can guarantee that Zeno behaviors
are avoided regardless of our choice of $\alpha$; however, we lose the
asymptotic convergence guarantee. That is, for $c_0 > 0$ we can only
guarantee convergence to a neighborhood of the desired average
consensus state.

As a result of the above discussion, we see that there is no way the
agents can choose the parameters $c_0, c_1,$ and $\alpha$ to ensure
asymptotic convergence to the average consensus state while also
guaranteeing Zeno executions are avoided. Consequently, more recent
works have proposed a local Lyapunov-based event-triggering condition
that only relies on currently available information and no exogenous
signals (e.g., time). This also naturally aligns when events are
triggered with the progression of the task as encoded in the Lyapunov
function. The state-dependent event-trigger to solve this problem was
first developed in~[20] and improved upon
in~[23] (removed global parameter~$a$ requirement); we
present this next.

Following~[23], to solve this problem we consider a
different Lyapunov function,
\begin{align*}
  V(x) = \frac{1}{2} (x - \bar{x} \mathbf{1})^T (x - \bar{x}
  \mathbf{1}) ,
\end{align*}
where $\bar{x} = \frac{1}{N} \sum_{i=1}^N x_i(0)$ is the average of
all initial conditions.  Then, given the closed-loop
dynamics~\eqref{eq:decentraldynamics}, we have
\begin{align*}
  \dot{V} = - x^T \dot{x} - \bar{x} \mathbf{1}^T \dot{x} = - x^T L
  \hat{x} - \bar{x} \mathbf{1}^T L \hat{x} = - x^T L \hat{x},
\end{align*}
where we have used the fact that the graph is weight-balanced in the
last equality.  As before, we are interested in finding conditions on
the error~$e$ such that $\dot{V} < 0$ at all times; however, we must
now do it without access to neighboring state information. Recalling
$e_i(t) = \hat{x}_i(t) - x_i(t)$, we can expand this out to
\begin{align*}
  \dot{V} &= -\hat{x}^T L \hat{x} + e^T L \hat{x} \\
  &= -\sum_{i=1}^N \sum_{j \in \NN_i} \left( \frac{1}{2} (\hat{x}_i -
    \hat{x}_j)^2 - e_i (\hat{x}_i-\hat{x}_j) \right) .
\end{align*}


Using Young's inequality for each product (see~[23]
for why this choice)
\begin{align*}
  e_i(\hat{x}_i - \hat{x}_j) \leq e_i^2 + \frac{1}{4} (\hat{x}_i -
  \hat{x}_j)^2
\end{align*}
yields
\begin{align*}
  \dot{V} &\leq -\sum_{i=1}^N \sum_{j \in \NN_i} \left( \frac{1}{2}
    (\hat{x}_i - \hat{x}_j)^2 - e_i^2 - \frac{1}{4} (\hat{x}_i -
    \hat{x}_j)^2 \right)
  \\
  &= - \sum_{i=1}^N \sum_{j \in \NN_i} \left(
    \frac{1}{4}(\hat{x}_i-\hat{x}_j)^2 - e_i^2 \right)
  \\
  &= \sum_{i=1}^N e_i | \NN_i | - \sum_{j \in \NN_i} \left(
    \frac{1}{4}(\hat{x}_i-\hat{x}_j)^2 \right) .
\end{align*}
If we can enforce the error of all agents to satisfy
\begin{align*}
  e_i^2 \leq \sigma_i \frac{1}{4 | \NN_i |} \sum_{j \in \NN_i}
  (\hat{x}_i - \hat{x}_j)^2
\end{align*}
with $\sigma_i \in (0,1)$ for all times, we have
\begin{align*}
  \dot{V} \leq \sum_{i=1}^N \frac{\sigma_i - 1}{4} \sum_{j \in \NN_i}
  (\hat{x}_i - \hat{x}_j)^2,
\end{align*}
which is strictly negative for all $L\hat{x} \neq 0$. The following
decentralized event-trigger ensures this is satisfied at all times.

\begin{theorem}[Decentralized event-triggered coordination
  (state-dependent)]\label{th:state}
  Given the closed-loop dynamics~\eqref{eq:decentraldynamics}, if the
  updates times of each agent~$i$ are determined as the times when
  \begin{align}
    f_i(e_i) \triangleq e_i^2 - \sigma_i \frac{1}{4 |\NN_i|} \sum_{j
      \in \NN_i} (\hat{x}_i - \hat{x}_j)^2 \geq 0,
  \end{align}
  then the system achieves multi-agent average consensus.
\end{theorem}

\begin{table}[htb]
  \centering
  \framebox[.9\linewidth]{\parbox{.85\linewidth}{%
      \parbox{\linewidth}{At times $t \in [t_\ell^i,t_{\ell+1}^i)$,
        agent~$i$ (continuously) performs:} 
      \vspace*{-1.5ex}
      \begin{algorithmic}[1]
      \STATE set $\hat{x}_i(t) = x_i(t_\ell^i)$
      \STATE set $e_i(t) = \hat{x}_i(t) - x_i(t)$
        \IF{$e_i(t)^2 \geq \sigma_i \frac{1}{4 | \NN_i |} \sum_{j \in \NN_i} (\hat{x}_i(t) - \hat{x}_j(t))^2$}  
        \STATE set $t_{\ell+1}^i = t$
        \STATE broadcast $\hat{x}_i(t) = x_i(t_{\ell+1}^i)$ to neighbors $j \in \NN_i$
        \STATE set $\ell = \ell + 1$
        \ENDIF
      \STATE set $u_i(t) = - \sum_{j \in \NN_i} (\hat{x}_i(t) - \hat{x}_j(t))$
      \end{algorithmic}}}
  \caption{Decentralized event-triggered coordination (state-dependent).}\label{tab:algorithm5}
\end{table}

It should be noted here that unlike all the other triggers presented
so far, this trigger is given by an inequality rather than an
equality. This is a result of the state-dependent triggering function
that agents use to determine when to communicate.  Since agents are
asynchronously sending each other messages, the information they have
about one another is also changing discontinuously.

\subsection{Directed graphs}\label{se:directed}

Up until now we have assumed that the communication graph was always
undirected. Here we extend the previous results to cases where the
communication graph~$\GG$ are directed but strongly connected and
weight-balanced.

More specifically, we say that agent~$i$ can only send messages to its
out-neighbors $j \in \NN_i^\text{out}$. Similarly, it can only receive
messages broadcast by its in-neighbors $j \in \NN_i^\text{in}$.
Conveniently, the closed-loop system dynamics is still given
by~\eqref{eq:decentraldynamics} where the only difference now is~$L$
is not symmetric. However, because it is weight-balanced we still have
that the sum of all states is an invariant quantity,
\begin{align*}
  \frac{d}{dt} \left( \mathbf{1}^T_N x(t) \right) = \mathbf{1}_N^T
  \dot{x}(t) = -\mathbf{1}_N^T L \hat{x}(t) = 0.
\end{align*}

\begin{remark}[Weight-balanced assumption] {\rm It should be noted
    that the weights of the directed graph for any digital
    implementations are design parameters that can be chosen to make a
    given directed communication topology weight-balanced. The
    works~[28,29] present provably correct
    distributed strategies that, given a directed communication
    topology, allow a network of agents to find such weight edge
    assignments. } \oprocend
\end{remark}

Remarkably, the same analysis from the previous section almost
directly follows and admits a similar triggering law. More
specifically, it can be shown that if we can enforce the error of all
agents to satisfy
\begin{align*}
  e_i^2 \leq \sigma_i \frac{1}{4 d_i^\text{out}} \sum_{j \in
    \NN_i^\text{out}} (\hat{x}_i - \hat{x}_j)^2 ,
\end{align*}
with $\sigma_i \in (0,1)$ for all times, we have
\begin{align}\label{eq:actual}
  \dot{V} \leq \sum_{i=1}^N \frac{\sigma_i - 1}{4} \sum_{j \in
    \NN_i^\text{out}} w_{ij} (\hat{x}_i - \hat{x}_j)^2,
\end{align}
which is strictly negative for all $L\hat{x} \neq 0$. The following
decentralized event-trigger ensures this is satisfied at all times.

\begin{theorem}[Decentralized event-triggered coordination
  on directed graphs]\label{th:state-directed}
  Given the closed-loop dynamics~\eqref{eq:decentraldynamics}, if the
  communication graph~$\GG$ is weight-balanced and the updates times
  of each agent~$i$ are determined as the times when
  \begin{align}
    f_i(e_i) \triangleq e_i^2 - \sigma_i \frac{1}{4 d_i^\text{out}}
    \sum_{j \in \NN_i} w_{ij} (\hat{x}_i - \hat{x}_j)^2 \geq 0,
  \end{align}
  then the system achieves multi-agent average consensus.
\end{theorem}

\begin{table}[htb]
  \centering
  \framebox[.9\linewidth]{\parbox{.85\linewidth}{%
      \parbox{\linewidth}{At times $t \in [t_\ell^i,t_{\ell+1}^i)$, agent~$i$ (continuously) performs:}
      \vspace*{-1.5ex}
      \begin{algorithmic}[1]
      \STATE set $\hat{x}_i(t) = x_i(t_\ell^i)$
      \STATE set $e_i(t) = \hat{x}_i(t) - x_i(t)$
        \IF{$e_i(t)^2 \geq \sigma_i \frac{1}{4 d_i^\text{out}} \sum_{j \in \NN_i^\text{out}} w_{ij}(\hat{x}_i(t) - \hat{x}_j(t))^2$}  
        \STATE set $t_{\ell+1}^i = t$
        \STATE broadcast $\hat{x}_i(t) = x_i(t_{\ell+1}^i)$ to in-neighbors $j \in \NN_i^\text{in}$
        \STATE set $\ell = \ell + 1$
        \ENDIF
      \STATE set $u_i(t) = - \sum_{j \in \NN_i^\text{out}} w_{ij} (\hat{x}_i(t) - \hat{x}_j(t))$
      \end{algorithmic}}}
  \caption{Decentralized event-triggered coordination on directed graphs.}\label{tab:algorithm6}
\end{table}

Unfortunately, most of the algorithms presented here are not
guaranteed to avoid Zeno behaviors making them risky to implement on
real systems. Moreover, the one algorithm that can in some cases
guarantee no Zeno behavior requires some global information.
In some cases modifications can be made to theoretically ensure
no Zeno behavior occurs; however, there may still be an arbitrarily
small amount of time between any two events (see e.g.,~[23])
making it undesirable from an implementation viewpoint. This is addressed in Remark~\ref{re:zeno} below
and the following section.

\begin{remark}[Implementation]\label{re:zeno} {\rm We note here an
    important issue regarding the connection between Zeno executions
    and implementation.  In general, dedicated hardware can only
    operate at some maximum frequency (e.g., a physical device can
    only broadcast a message or evaluate a function a finite number of
    times in any finite period of time). This means that ensuring a
    system does not exhibit Zeno behavior may not be enough to
    guarantee the algorithm can be implemented on a physical system if
    the physical hardware cannot match the speed of actions required
    by the algorithm. More specifically, it is guaranteed that Zeno
    behavior does not exist if the sequence of times $t_\ell^i
    \rightarrow \infty$ as $\ell \rightarrow \infty$; however, this is
    not as strong as ensuring that there exists a minimum time in
    between triggers $t_{\ell+1}^i - t_\ell^i \geq \tau^\text{min} >
    0$, which is a more pragmatic constraint when considering physical
    hardware.  } \oprocend
\end{remark}

In light of Remark~\ref{re:zeno}, we consider enforcing a minimum time
between events in the following section.

\section{Periodic event-triggered coordination}

Throughout this chapter we have assumed that all event-triggers can be
evaluated continuously.  That is, the exact moment at which a
triggering condition is met, an action (e.g., state broadcast and
control signal update) is carried out. However, this may still be an
unrealistic assumption when considering digital implementations. More
specifically, a physical device cannot continuously evaluate whether a
triggering condition has occurred or not. This observation motivates
the need for studying sampled-data (or periodically checked)
event-triggered coordination strategies.

Specifically, given a sampling period $h \in \realpositive$, we let
$\{ t_{\ell'} \}_{\ell' \in \integernonnegative}$, where $t_{\ell'+1}
= t_\ell' + h$, denote the sequence of times at which agents evaluate
the decision of whether to broadcast their state to their
neighbors. This type of design is more in line with the constraints
imposed by real-time implementations, where individual components work
at some fixed frequency, rather than continuously. An inherent and
convenient feature of this strategy is the lack of Zeno behavior
(since inter-event times are naturally lower bounded by~$h$).

Consequently, we begin by revisiting the result of
Theorem~\ref{th:state-directed}. Intuitively, as long as the sampling
period~$h$ is small enough, the closed-loop system with a periodically
checked event-triggering condition will behave similarly to the system
with triggers being evaluated continuously. The proof of convergence
for the triggering law in Theorem~\ref{th:state-directed} hinges on
the fact that
\begin{align*}
  e_i^2(t) \leq \sigma_i \frac{1}{4 d_i^\text{out}} \sum_{j \in \NN_i}
  w_{ij} (\hat{x}_i(t) - \hat{x}_j(t))^2
\end{align*}
for all times~$t$. Instead, since we now assume the triggering
function~\ref{eq:decentraltrigger} can only be evaluated periodically,
we have that
\begin{align}\label{eq:enforceperiodic}
  e_i^2(t_{\ell'}) \leq \sigma_i \frac{1}{4 d_i^\text{out}} \sum_{j
    \in \NN_i} w_{ij} (\hat{x}_i(t_{\ell'}) - \hat{x}_j(t_{\ell'}))^2
\end{align}
is only guaranteed at the specific times $\{ t_{\ell'} \}_{\ell' \in
  \integernonnegative}$ at which the triggering function can be
evaluated. The algorithm is formalized in Table~\ref{tab:algorithm7}.

It should be noted that this algorithm is identical to the one in
Table~\ref{tab:algorithm6} except it is only executed periodically now
rather than continuously.  The following result then provides a
sufficient condition on how small the period~$h$ has to be to
guarantee convergence. The result is obtained by analyzing what
happens to the Lyapunov function~$V$ in between these
times. 

\begin{theorem}[Periodic event-triggered coordination]\label{th:sampled-data}
  Given the closed-loop dynamics~\eqref{eq:decentraldynamics}, if the
  communication graph~$\GG$ is weight-balanced and the update times of
  each agent~$i$ are determined as the times $t' \in \{0, h, 2h, \dots
  \}$ when
  \begin{align*}
  f_i(e_i) \triangleq e_i^2 - \sigma_i \frac{1}{4 d_i^\text{out}} \sum_{j \in
      \NN_i} w_{ij} (\hat{x}_i - \hat{x}_j)^2 \geq 0,
  \end{align*}
  and $h \in \realpositive$ satisfies
  \begin{align}\label{eq:bound}
    \sigma_{\max} + 4 h w_{\max} |\Noutmax | < 1,
  \end{align}
  where $ w_{\max} = \max_{i \in \until{N}, j \in \NN_i^\text{out}} w_{ij}$ and
  $|\Noutmax | = \max_{i \in \until{N}} |\Nouti |$,
  then the system achieves multi-agent average consensus.
\end{theorem}

Note that checking the sufficient condition~\eqref{eq:bound} requires
knowledge of the global quantities $\sigma_{\max}$, $w_{\max}$,
and~$\Noutmax$. Ensuring that this condition is met can either be
enforced a priori by the designer or, alternatively, the network can
execute a distributed initialization procedure,
e.g.,~[30,3], to compute these quantities in finite
time.  Once known, agents can compute $h$ by instantiating a specific
formula to select it that is guaranteed to satisfy~\eqref{eq:bound}.

\begin{table}[htb]
  \centering
  \framebox[.9\linewidth]{\parbox{.85\linewidth}{%
      \parbox{\linewidth}{At times $t \in \{0, h, 2h, \dots \}$, agent~$i$ performs:}
      \vspace*{-1.5ex}
      \begin{algorithmic}[1]
      \STATE set $\hat{x}_i(t) = x_i(t_\ell^i)$
      \STATE set $e_i(t) = \hat{x}_i(t) - x_i(t)$
        \IF{$e_i(t)^2 \geq \sigma_i \frac{1}{4 d_i^\text{out}} \sum_{j \in \NN_i^\text{out}} w_{ij}(\hat{x}_i(t) - \hat{x}_j(t))^2$}  
        \STATE set $t_{\ell+1}^i = t$
        \STATE broadcast $\hat{x}_i(t) = x_i(t_{\ell+1}^i)$ to neighbors $j \in \NN_i^\text{in}$
        \STATE set $\ell = \ell + 1$
        \ENDIF
      \STATE set $u_i(t) = - \sum_{j \in \NN_i^\text{out}} w_{ij} (\hat{x}_i(t) - \hat{x}_j(t))$
      \end{algorithmic}}}
  \caption{Periodic event-triggered coordination on directed graphs.}\label{tab:algorithm7}
\end{table}

A drawback of this algorithm is that the period~$h$ must be the same
for all agents, requiring synchronous action. It is not difficult to
envision asynchronous versions of this algorithm for which correctness
guarantees have not currently been established.

\section{Conclusions and future outlook}

This chapter has presented a high-level overview of the ideas behind
event-triggered communication and control applied to multi-agent
average consensus problems.  Although Table~\ref{ta:survey} makes it
look like there is a complete story concerning event-triggered
consensus problems, this is certainly not true as there still remain
many issues to be addressed regarding asynchronism, guarantees on
non-Zeno behavior, and practical considerations. There are indeed
still exciting new directions being explored at the time of writing
that would only serve to expand this table in the future.  For
instance, all event-triggered protocols discussed in this chapter
assume that all agents are able to `listen' for incoming messages at
all times. In other words, when a message is broadcast by an
agent~$i$, this message is immediately received by a neighboring
agent~$j \in \NN_i$ who immediately (or within some reasonable time
due to delays, etc.) reacts to this event by updating its control
signal.  However, this may not be possible in all scenarios which
presents a whole new set of technical challenges. For example, some
recent preliminary results have been developed in this setup motivated
by the need for coordinating
submarines~[31,32], where agents can only
communicate when they are at the surface of the water. While a
submarine is submerged, any message broadcast by another submarine
cannot be received until it resurfaces.


\section*{References}
\begin{enumerate}
\item 
R.~Olfati-Saber, J.~A. Fax, and R.~M. Murray, ``Consensus and cooperation in
  networked multi-agent systems,'' \emph{Proceedings of the IEEE}, vol.~95,
  no.~1, pp. 215--233, 2007.

\item 
F.~Bullo, J.~Cort{\'e}s, and S.~Mart{\'\i}nez, \emph{Distributed Control of
  Robotic Networks}, ser. Applied Mathematics Series.\hskip 1em plus 0.5em
  minus 0.4em\relax Princeton University Press, 2009, electronically available
  at \url{http://coordinationbook.info}.

\item
W.~Ren and R.~W. Beard, \emph{Distributed Consensus in Multi-Vehicle
  Cooperative Control}, ser. Communications and Control Engineering.\hskip 1em
  plus 0.5em minus 0.4em\relax Springer, 2008.

\item
M.~Mesbahi and M.~Egerstedt, \emph{Graph Theoretic Methods in Multiagent
  Networks}, ser. Applied Mathematics Series.\hskip 1em plus 0.5em minus
  0.4em\relax Princeton University Press, 2010.

\item
K.~J. {\r Astr{\"o}m} and B.~M. Bernhardsson., ``Comparison of {R}iemann and
  {L}ebesgue sampling for first order stochastic systems,'' in \emph{{IEEE}
  Conf.\ on Decision and Control}, Las Vegas, NV, Dec. 2002, pp. 2011--2016.

\item
W.~P. M.~H. Heemels, K.~H. Johansson, and P.~Tabuada, ``An introduction to
  event-triggered and self-triggered control,'' in \emph{{IEEE} Conf.\ on
  Decision and Control}, Maui, HI, 2012, pp. 3270--3285.

\item
M.~{Mazo Jr.} and P.~Tabuada, ``Decentralized event-triggered control over
  wireless sensor/actuator networks,'' \emph{IEEE Transactions on Automatic
  Control}, vol.~56, no.~10, pp. 2456--2461, 2011.

\item
X.~Wang and M.~D. Lemmon, ``Event-triggering in distributed networked control
  systems,'' \emph{IEEE Transactions on Automatic Control}, vol.~56, no.~3, pp.
  586--601, 2011.

\item
C.~St\"oker, D.~Vey, and J.~Lunze, ``Decentralized event-based control:
  Stability analysis and experimental evaluation,'' \emph{Nonlinear Analysis:
  Hybrid Systems}, vol.~10, pp. 141--155, 2013.

\item
G.~Xie, H.~Liu, L.~Wang, and Y.~Jia, ``Consensus in networked multi-agent
  systems via sampled control: fixed topology case,'' in \emph{{A}merican
  {C}ontrol {C}onference}, St. Louis, MO, 2009, pp. 3902--3907.

\item
W.~P. M.~H. Heemels and M.~C.~F. Donkers, ``Model-based periodic
  event-triggered control for linear systems,'' \emph{Automatica}, vol.~49,
  no.~3, pp. 698--711, 2013.

\item
X.~Meng and T.~Chen, ``Event based agreement protocols for multi-agent
  networks,'' \emph{Automatica}, vol.~49, no.~7, pp. 2125--2132, 2013.

\item
M.~Zhong and C.~G. Cassandras, ``Asynchronous distributed optimization with
  event-driven communication,'' \emph{IEEE Transactions on Automatic Control},
  vol.~55, no.~12, pp. 2735--2750, 2010.

\item
A.~Anta and P.~Tabuada, ``To sample or not to sample: self-triggered control
  for nonlinear systems,'' \emph{IEEE Transactions on Automatic Control},
  vol.~55, no.~9, pp. 2030--2042, 2010.

\item
X.~Wang and M.~D. Lemmon, ``Self-triggered feedback control systems with
  finite-gain {L}$_2$ stability,'' \emph{IEEE Transactions on Automatic
  Control}, vol.~54, no.~3, pp. 452--467, 2009.

\item
C.~Nowzari and J.~Cort{\'e}s, ``Team-triggered coordination for real-time
  control of networked cyberphysical systems,'' \emph{IEEE Transactions on
  Automatic Control}, vol.~61, no.~1, pp. 34--47, 2016.

\item
D.~V. Dimarogonas, E.~Frazzoli, and K.~H. Johansson, ``Distributed
  event-triggered control for multi-agent systems,'' \emph{IEEE Transactions on
  Automatic Control}, vol.~57, no.~5, pp. 1291--1297, 2012.

\item
Z.~Liu, Z.~Chen, and Z.~Yuan, ``Event-triggered average-consensus of
  multi-agent systems with weighted and direct topology,'' \emph{Journal of
  Systems Science and Complexity}, vol.~25, no.~5, pp. 845--855, 2012.

\item
G.~S. Seyboth, D.~V. Dimarogonas, and K.~H. Johansson, ``Event-based
  broadcasting for multi-agent average consensus,'' \emph{Automatica}, vol.~49,
  no.~1, pp. 245--252, 2013.

\item
E.~Garcia, Y.~Cao, H.~Yu, P.~Antsaklis, and D.~Casbeer, ``Decentralised
  event-triggered cooperative control with limited communication,''
  \emph{International Journal of Control}, vol.~86, no.~9, pp. 1479--1488,
  2013.

\item
G.~Guo, L.~Ding, and Q.~Han, ``A distributed event-triggered transmission
  strategy for sampled-data consensus of multi-agent systems,''
  \emph{Automatica}, vol.~50, pp. 1489--1496, 2014.

\item
X.~Meng, L.~Xie, Y.~C. Soh, C.~Nowzari, and G.~J. Pappas, ``Periodic
  event-triggered average consensus over directed graphs,'' in \emph{{IEEE}
  Conf.\ on Decision and Control}, Osaka, Japan, Dec. 2015, pp. 4151--4156.

\item
C.~Nowzari and J.~Cort\'es, ``Distributed event-triggered coordination for
  average consensus on weight-balanced digraphs,'' \emph{Automatica}, vol.~68,
  pp. 237--244, 2016.

\item
M.~{Mazo Jr.}, A.~Anta, and P.~Tabuada, ``On self-triggered control for linear
  systems: Guarantees and complexity,'' in \emph{European Control Conference},
  Budapest, Hungary, Aug. 2009, pp. 3767--3772.

\item
R.~Olfati-Saber and R.~M. Murray, ``Consensus problems in networks of agents
  with switching topology and time-delays,'' \emph{IEEE Transactions on
  Automatic Control}, vol.~49, no.~9, pp. 1520--1533, 2004.

\item
R.~Aragues, G.~Shi, D.~V. Dimaragonas, C.~Sagues, and K.~H. Johansson,
  ``Distributed algebraic connectivity estimation for adaptive event-triggered
  consensus,'' in \emph{{A}merican {C}ontrol {C}onference}, Montreal, Canada,
  2012, pp. 32--37.

\item
P.~Yang, R.~A. Freeman, G.~J. Gordon, K.~M. Lynch, S.~S. Srinivasa, and
  R.~Sukthankar, ``Decentralized estimation and control of graph connectivity
  for mobile sensor networks,'' \emph{Automatica}, vol.~46, no.~2, pp.
  390--396, 2010.

\item
B.~Gharesifard and J.~Cort{\'e}s, ``Distributed strategies for generating
  weight-balanced and doubly stochastic digraphs,'' \emph{European Journal of
  Control}, vol.~18, no.~6, pp. 539--557, 2012.

\item
A.~Rikos, T.~Charalambous, and C.~N. Hadjicostis, ``Distributed weight
  balancing over digraphs,'' \emph{{IEEE} Transactions on Control of Network
  Systems}, 2014, to appear.

\item
N.~A. Lynch, \emph{Distributed Algorithms}.\hskip 1em plus 0.5em minus
  0.4em\relax Morgan Kaufmann, 1997.

\item
A.~Adaldo, D.~Liuzza, D.~V. Dimarogonas, and K.~H. Johansson, ``Control of
  multi-agent systems with event-triggered cloud access,'' in \emph{{E}uropean
  {C}ontrol {C}onference}, Linz, Austria, 2015, pp. 954--961.

\item
C.~Nowzari and G.~J. Pappas, ``Multi-agent coordination with asynchronous cloud
  access,'' in \emph{{A}merican {C}ontrol {C}onference}, Boston, MA, Jun. 2016,
  to appear.

\item
H.~K. Khalil, \emph{Nonlinear Systems}, 3rd~ed.\hskip 1em plus 0.5em minus
  0.4em\relax Prentice Hall, 2002.

\end{enumerate}
\backmatter

\appendix

\section*{Appendix}

\paragraph{\textbf{Proof of Theorem~\ref{th:main}}}
Consider the Lyapunov function
\begin{align*}
  V(x) = \frac{1}{2} x^T L x .
\end{align*}
Then, given the dynamics~\eqref{eq:dynamics} and the continuous
control law~\eqref{eq:ideal},
\begin{align*}
  \dot{V}(x) = x^T L \dot{x} = -x^T L^T L x = - \TwoNorm{Lx}^2,
\end{align*}
where we have used the fact that~$L$ is symmetric. It is now clear
that using the continuous control law~\eqref{eq:ideal} we have
$\dot{V}(x) < 0$ for all $Lx \neq 0$. Using LaSalle's Invariance
Principle~[33], it can then be shown that
\begin{align*}
  x(t) \rightarrow \{ Lx = 0 \} = \{ x_i = x_j \forall i,j \in \until{N} \}
\end{align*}
as $t \rightarrow \infty$. Combining this with the fact that the sum
of all states is an invariant quantity concludes the proof,
\begin{align*}
  \frac{d}{dt} \left( \ones{N}^T x(t) \right) = \ones{N}^T \dot{x}(t)
  = -\ones{N}^T L x(t) = 0 .
\end{align*}
\hfill $\blacksquare$

\paragraph{\textbf{Proof of Theorem~\ref{th:time}}}
Let $\delta(t) = x(t) - \bar{x} \mathbf{1}$, where $\bar{x} = \frac{1}{N}
\sum_{i=1}^N x_i(0)$ is the average of all initial conditions. Then,
$\dot{\delta}(t) = - L \delta(t) - L e(t)$, yielding
\begin{align*}
  \delta(t) = e^{-Lt} \delta(0) - \int_0^t e^{-L(t-s)} L e(s) ds .
\end{align*}
Taking norms,
\begin{align*}
  \TwoNorm{\delta(t)} &\leq \TwoNorm{ \delta(0) e^{-Lt}} + \int_0^t
  \TwoNorm{ e^{-L(t-s)} Le(s) } ds
  \\
  &\leq e^{- \lambda_2(L) t} \TwoNorm{\delta(0)} + \int_0^t
  e^{-\lambda_2(L)(t-s)} \TwoNorm{L e(s)} ds, 
\end{align*}
where the second inequality follows from~[19, Lemma 2.1].

Using the condition
\begin{align*}
  | e_i(t) | \leq c_0 + c_1 e^{-\alpha t},
\end{align*}
it follows that
\begin{align*}
  \TwoNorm{\delta(t)} &\leq e^{-\lambda_2 t} \TwoNorm{\delta(0)} +
  \TwoNorm{L} \sqrt{N} \int_0^t e^{-\lambda_2(t-s)} (c_0 + c_1
  e^{-\alpha s}) ds
  \\
  &= e^{-\lambda_2 t} \left( \TwoNorm{\delta(0)} - \TwoNorm{L}
    \sqrt{N} \left( \frac{c_0}{\lambda_2} + \frac{c_1}{\lambda_2 -
        \alpha} \right) \right) + e^{-\alpha t} \frac{ \TwoNorm{L}
    \sqrt{N} c_1}{\lambda_2 - \alpha} + \frac{ \TwoNorm{L} \sqrt{N}
    c_0 }{\lambda_2}.
\end{align*}
The convergence result then follows by taking $t \rightarrow \infty$. 
\hfill $\blacksquare$

\paragraph{\textbf{Proof of Theorem~\ref{th:sampled-data}}}
Since~\eqref{eq:enforceperiodic} is only guaranteed at the sampling
times under the periodic event-triggered coordination algorithm presented
in Table~\ref{tab:algorithm7}, we analyze what happens to the Lyapunov
function~$V$ in between them.  For $t \in [t_{\ell'}, t_{\ell'+1})$, note
that
\begin{align*}
  e(t) = e(t_{\ell'}) + (t-t_{\ell'}) L \hat{x}(t_{\ell'}) .
\end{align*}
Substituting this expression into $\dot{V}(t) = -\hat{x}^T(t) L
\hat{x}(t) + e^T(t) L \hat{x}(t)$, we obtain
\begin{multline*}
  \dot{V}(t) = -\hat{x}^T(t_{\ell'}) L \hat{x}(t_{\ell'}) + e^T(t_{\ell'}) L
  \hat{x}(t_{\ell'})
  \\
  + (t - t_{\ell'}) \hat{x}^T(t_{\ell'}) L^T L \hat{x}(t_{\ell'}) ,
\end{multline*}
for all $t \in [t_{\ell'}, t_{\ell'+1})$.  For a simpler exposition,
we drop all arguments referring to time $t_{\ell'}$ in the sequel.
Following a similar discussion to Section~\ref{se:coordination}, it can be
shown that
\begin{align*}
  \dot{V}(t) \leq \sum_{i=1}^N \frac{\sigma_i - 1}{4} \sum_{j \in
    \NN_i^\text{out}} w_{ij} (\hat{x}_i - \hat{x}_j)^2 + (t-t_{\ell'})
  \hat{x}^T L^T L \hat{x} .
\end{align*}
Note that the first term is exactly what we have when we are able to
monitor the triggers continuously~\eqref{eq:actual}.

Using the fact that $\left( \sum_{k=1}^p y_k \right)^2 \leq p
\sum_{k=1}^p y_k^2$ (which follows directly from the Cauchy-Schwarz
inequality), we bound
\begin{align}\label{eq:boundxx}
  \hat{x}^T L^T L \hat{x} & = \sum_{i=1}^N \left( \sum_{j \in \Nouti}
  w_{ij} (\hat{x}_i-\hat{x}_j) \right)^2 \notag
  \\
  & \le \sum_{i=1}^N |\Nouti| w^{\max}_i \sum_{j \in \Nouti} w_{ij} (\hat{x}_i -
  \hat{x}_j)^2 \notag
  \\
  & = |\Noutmax | w_{\max} \sum_{i=1}^N \sum_{j \in \Nouti} w_{ij} (\hat{x}_i -
  \hat{x}_j)^2 ,
\end{align}
where $w^{\max}_i = \max_{j \in \NN_i^\text{out}} w_{ij}$. 
Hence, for $t \in [t_{\ell'}, t_{\ell'+1})$,
\begin{align*}
  \dot{V}(t) & \leq \sum_{i=1}^N \Big( \frac{\sigma_i - 1}{4} + h
  w_{\max} |\Noutmax | \Big) \sum_{j \in \Nouti} w_{ij} (\hat{x}_i -
  \hat{x}_j)^2
  \\
  & \leq \Big( \frac{\sigma_{\max}-1}{2} + 2 h
  w_{\max} |\Noutmax | \Big) \hat{x}^T L \hat{x} .
\end{align*}
Then, by using~\eqref{eq:bound}, it can be shown that there exists $\BB > 0$
such that
\begin{align*}
  \dot{V}(t) & \leq \frac{1}{2 \BB }\Big( \sigma_{\max} + 4 h w_{\max}
  |\Noutmax | - 1 \Big) V(x(t)),
\end{align*}
which implies the result. See~[23] for details.
\hfill $\blacksquare$

\end{document}